\documentstyle[art11]{article}

\newdimen\mypt
\newdimen\halfpt
\def\pic#1 #2 #3;{\raisebox{0.6ex}{\raisebox{-#2\halfpt}{\vbox
      to #2\mypt{\hbox to #3\mypt{\special{em:graph
      #1.bmp}\hfill}\vfill}}}}

\def\proved{\ifmmode\eqno\Box\medskip\else
        \nobreak\hfill\nopagebreak\discretionary{}
        {\hbox to\textwidth{\hfill$\Box$}}{\hbox{$\Box$}}\par
        \addvspace\medskipamount\fi}

\newcounter{defin}

\newcounter{remark}

\font\msamfive=msam5 at 3pt
\def\uparc(#1_#2#3_#4){{\mathop{#1_{#2}#3_{#4}}\limits^{\displaystyle
\frown\kern-3pt\llap{\raisebox{1.5pt}{\msamfive\char121}}
\phantom{\scriptscriptstyle#4}
}}}
\def\loarc(#1_#2#3_#4){{\mathop{#1_{#2}#3_{#4}}\limits_{\displaystyle
\smile\kern-3pt\llap{\raisebox{1.7pt}{\msamfive\char113}}
\phantom{\scriptscriptstyle#4}
}}}

\def\sing{{\rm sing \,}}

\author{S. D. Tyurina\thanks{This work is partly
supported by INTAS, grant no.~97-1644.}
}
\date{Dep. of Phys. and Math.\\
Kolomensky Pedagogical Institute\\
Kolomna 140411 RUSSIA\\
e-mail: tyurina@mccme.ru}
\title{{\sc Diagram formulas of the Viro-Polyak type for finit degree
invariants.}\\
{\scriptsize(The brief version of this paper is in UMN {\bf t.54} N3 (1998).
)}}

\begin{document}
\baselineskip=11.8pt
\maketitle
UDK 515.162
\abstract{
Vassiliev's knot invariants can be computed in different ways but many of them
as Kontsevich integral are very difficult. We consider more visual diagram
formulas of the type Polyak-Viro and give new diagram formula for the two
basic Vassiliev invariant of degree 4.}
\section{The Vassiliev knot invariants.}
Let $K:S^1 \longrightarrow {\bf R^3}$
be an oriented knot and $K^{\sing}_n:S^1
\longrightarrow {\bf R^3}$ be a singular knot
with $n$ double points.
Denote by ${\cal K}$ the space of knots and
by ${\cal K}^{\sing}$ the space of singular knots.
Any knot invariant $V:{\cal K} \longrightarrow {\bf Q}$ may be extended from
ordinary knots to singular knots by next inductive rule:

\begin{picture}(420,22)
\put(80,9){$V^{(i)}($}
\put(120,11){\circle*{2}}
\put(132,9){$)=V^{(i-1)}($}
\put(217,9){$)-V^{(i-1)}($}
\put(299,9){).}
\put(420,9){(1)}
\put(109,0){\vector(1,1){22}}
\put(131,0){\vector(-1,1){22}}
\put(193,0){\vector(1,1){22}}
\put(215,0){\line(-1,1){10}}
\put(203,12){\vector(-1,1){10}}
\put(274,0){\line(1,1){10}}
\put(288,12){\vector(1,1){10}}
\put(297,0){\vector(-1,1){22}}
\end{picture}

{\bf Definition.} A knot invariant $V:{\cal K} \longrightarrow
{\bf Q}$ is said to be Vassiliev invariant of degree
less then or equal to $n$, if there exists $n \in {\bf N}$ such that
$$V^{(n+1)} \equiv 0.$$

\section{Chord diagrams.}
A chord diagram $D_n$ of degree $n$ is the circle with a distinguished set
of $n$ unordered pairs of points regarded up to orientation preserving
diffeomorphisms of the circle. Denote by ${\cal D}_n$ the space generated
by chord diagrams of degree n over {\bf Q}.
\section{Weight systems.}
{\bf Definition.} A linear function $W$ is called a weight system of degree $n$
if it satisfies 1-term and 4-term relations.
For any invariant $f$ we construct homomorphism $W_f:{\cal D}_n
\longrightarrow {\bf Q}$ as follows:
$$W_f = f|_{{\cal K}^{sing}_n}.$$
The homomorphism $W_f$ is the weight system.

\section{The Gauss diagrams.}
A chord diagram of the singular knot is the circle with pre-images
of double points connected with  chords.
To obtain the analogus diagram
of an ordinary knot (that is called the Gauss diagram) from the chord diagram
of corresponding singular knot we must give the information on overpasses and
underpasses. Each chord is oriented from the upper branch to the lower one
and equipped with the sign (the local writh number of corresponding double
point of planar projection of the knot).

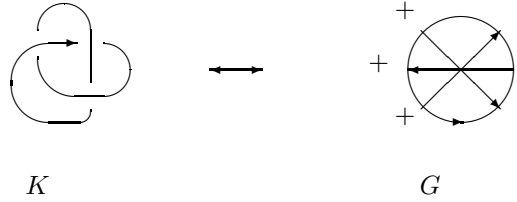
\begin{picture}(400,65)
\put(100,45){\vector(1,0){5}}
\put(100,50){\oval(20,20)[t]}
\put(110,30){\line(0,1){20}}
\put(105,30){\oval(50,30)[tl]}
\put(110,40){\oval(40,30)[bl]}
\put(110,25){\line(1,0){5}}
\put(95,20){\oval(30,10)[br]}
\put(115,35){\oval(20,20)[r]}
\put(95,30){\oval(30,30)[bl]}

\put(165,35){\vector(-1,0){10}}
\put(165,35){\vector(1,0){10}}
\put(250,35){\circle{40}}
\put(250,15){\vector(1,0){1}}
\put(270,35){\vector(-1,0){40}}
\put(235,20){\vector(1,1){30}}
\put(235,50){\vector(1,-1){30}}
\put(215,35){$+$}
\put(225,15){$+$}
\put(225,55){$+$}
\end{picture}

\hskip3cm {$K$\qquad \qquad \qquad \qquad \qquad \qquad \qquad $G$}

\section{Formula of the type of Viro-Polyak for Vassiliev invariants of degree 4.}
Denote by $<A,G>$ algebraic number of subdiagrams of given combinatorial
type $A \quad A \subset G$ and let $<\sum_{i}n_iA_i,G> = \sum_{i}n_i<A_i,G>,
\quad n_i \in {\bf Q}$ by definition.

\noindent
{\bf Theorem.} {\it Let $W_4$ be a weight system such that $W_4(
\begin{picture}(20,15)
\put(8,3){\circle{15}}     
\put(3,-2){\line(1,1){10}}  
\put(13,-2){\line(-1,1){10}} 
\put(8,-5){\line(0,1){15}}  
\put(1,3){\line(1,0){15}}  
\end{picture}
)=0.$
Then following formula gives Vassiliev invariant of degree 4:}

\begin{picture}(420,40)
\put(4,10){\it $f_4(K)=\frac{1}{2}W_4($}
\put(100,15){\circle{30}}      
\put(95,1){\line(-1,2){10}}   
\put(86,7){\line(1,2){12}}     
\put(105,1){\line(1,2){10}}    
\put(114,7){\line(-1,2){12}}   
\put(120,10){$)<$}
\put(160,15){\circle{30}}      
\put(146,23){\vector(1,-2){12}}   
\put(146,7){\vector(1,2){12}}     
\put(165,1){\vector(1,2){10}}    
\put(163,31){\vector(1,-2){12}}   
\put(180,10){$,G>+$}

\end{picture}

\begin{picture}(420,40)
\put(52,10){$\frac{1}{4}W_4($}
\put(100,15){\circle{30}}      
\put(90,2){\line(0,1){25}}   
\put(100,-1){\line(0,1){32}}   
\put(84,15){\line(1,0){32}}   
\put(110,2){\line(0,1){25}}   
\put(120,10){$)<$}

\put(160,15){\circle{30}}      
\put(150,28){\vector(0,-1){26}}   
\put(170,28){\vector(0,-1){26}}   
\put(144,15){\vector(1,0){32}}   
\put(160,-1){\vector(0,1){32}}   
\put(180,10){+}
\put(210,15){\circle{30}}        
\put(200,2){\vector(0,1){26}}   
\put(220,2){\vector(0,1){26}}   
\put(194,15){\vector(1,0){32}}   
\put(210,31){\vector(0,-1){32}} 
\put(230,10){$,G>+$}
\end{picture}

\begin{picture}(420,40)
\put(52,10){$\frac{1}{2}W_4($}
\put(100,15){\circle{30}}     
\put(120,10){$)<$}           
\put(93,1){\line(0,1){28}}  
\put(107,1){\line(0,1){28}}  
\put(86,7){\line(1,0){28}}  
\put(86,23){\line(1,0){28}} 

\put(160,15){\circle{30}}      
\put(153,30){\vector(0,-1){29}}  
\put(167,1){\vector(0,1){28}}  
\put(146,7){\vector(1,0){28}}  
\put(174,23){\vector(-1,0){28}} 
\put(180,10){+}
\put(210,15){\circle{30}}        
\put(230,10){$,G>+$}           
\put(203,1){\vector(0,1){28}}  
\put(217,29){\vector(0,-1){28}}  
\put(224,7){\vector(-1,0){28}}  
\put(196,23){\vector(1,0){28}} 
\end{picture}

\begin{picture}(420,40)
\put(52,10){$\frac{1}{4}W_4($}           
\put(100,15){\circle{30}}     
\put(120,10){$)<$}           
\put(100,-1){\line(-1,2){13}}  
\put(87,5){\line(1,2){13}}   
\put(93,0){\line(2,1){22}}  
\put(93,30){\line(2,-1){22}}

\put(160,15){\circle{30}}      
\put(180,10){+}                   
\put(160,-1){\vector(-1,2){13}}    
\put(160,31){\vector(-1,-2){13}}  
\put(175,11){\vector(-2,-1){22}}  
\put(154,30){\vector(2,-1){22}}  
\put(210,15){\circle{30}}        

\put(230,10){+}                   
\put(197,25){\vector(1,-2){13}}    
\put(197,5){\vector(1,2){13}}  
\put(225,12){\vector(-2,-1){22}}  
\put(203,29){\vector(2,-1){22}}  

\put(260,15){\circle{30}}         
\put(330,10){,$G>+$}                   
\put(248,25){\vector(1,-2){13}}    
\put(250,3){\vector(1,2){14}}  
\put(254,1){\vector(2,1){22}}  
\put(276,18){\vector(-2,1){22}}  

\put(310,15){\circle{30}}         
\put(280,10){+}                   
\put(310,-1){\vector(-1,2){13}}    
\put(310,31){\vector(-1,-2){13}}  
\put(303,1){\vector(2,1){22}}  
\put(325,19){\vector(-2,1){22}}  
\end{picture}

\begin{picture}(420,40)
\put(52,10){$\frac{1}{4}W_4($}           
\put(100,15){\circle{30}}     
\put(120,10){$)<$}           
\put(109,1){\line(0,1){27}}  
\put(84,15){\line(1,0){31}}  
\put(100,-1){\line(-1,2){13}}  
\put(87,5){\line(1,2){13}}   

\put(160,15){\circle{30}}      
\put(169,1){\vector(0,1){27}}  
\put(144,15){\vector(1,0){31}}  
\put(180,10){+}                   
\put(160,-1){\vector(-1,2){13}}    
\put(160,31){\vector(-1,-2){13}}  

\put(210,15){\circle{30}}        
\put(219,28){\vector(0,-1){27}}  
\put(225,15){\vector(-1,0){31}}  
\put(230,10){+}                   
\put(197,25){\vector(1,-2){13}}    
\put(197,5){\vector(1,2){13}}  

\put(260,15){\circle{30}}         
\put(269,1){\vector(0,1){27}}  
\put(275,15){\vector(-1,0){31}}  
\put(330,10){,$G>+$}                   
\put(248,25){\vector(1,-2){13}}    
\put(248,5){\vector(1,2){13}}  

\put(310,15){\circle{30}}         
\put(319,28){\vector(0,-1){27}}  
\put(294,15){\vector(1,0){31}}  
\put(280,10){+}                   
\put(310,-1){\vector(-1,2){13}}    
\put(310,31){\vector(-1,-2){13}}  
\put(330,10){,$G>+$} 
\end{picture}

\begin{picture}(420,40)
\put(52,10){$\frac{1}{4}W_4($} 
\put(100,15){\circle{30}}     
\put(120,10){$)<$}           
\put(94,1){\line(0,1){29}}  
\put(106,1){\line(0,1){29}}  
\put(86,8){\line(2,1){29}}  
\put(86,22){\line(2,-1){29}}

\put(160,15){\circle{30}}      
\put(180,10){+}                   
\put(154,1){\vector(0,1){29}}  
\put(166,30){\vector(0,-1){29}}  
\put(146,8){\vector(2,1){29}}  
\put(174,8){\vector(-2,1){29}}

\put(210,15){\circle{30}}        
\put(204,30){\vector(0,-1){29}}  
\put(216,1){\vector(0,1){29}}  
\put(196,8){\vector(2,1){29}}  
\put(224,8){\vector(-2,1){29}}
\put(230,10){,$G>,$}
\end{picture}

\noindent
{\it where  $G$ is a Gauss diagram of a knot $K$.}

\noindent
{\bf Proof} contains the checking of independens of formula under the
Reidemeister moves.

Denote by $V_4^1$ and $V_4^2$ the Vassiliev invariants corresponding to basic
weight systems satisfying theorem. $V_4^1$, $V_4^2$ and $V_4$ from [2] (see thm.4) is
basis in the space of the Vassiliev invariants of degree 4.

\vskip1cm
Svetlana D. Tyurina

Kolomensky Pedagogical Institute

140411, Moscow reg., Kolomna, Zelenaja 30

e-mail: ktti@kolomna.ru

\begin{thebibliography}{99}
\bibitem{Lannes}
J.Lannes.
Sur les invariants de Vassiliev de degre$\acute e$ inferieur
ou $\acute e$gal $\grave a$ 3.
{\it L'En\-seignement Math$\acute e$matique} {\bf t.39} (1993), pp.~295--316.
\bibitem{VP}
M.Polyak, O.Viro.
Gauss diagram formulae for Vassiliev invariants.
{\it Int. Math. Res. Notices} {\bf 11} (1994), pp.~445-453.
\bibitem{Vas}
V.A.Vassiliev.
Complemets of discriminants of smoos maps.
{\it Amer. Math. Soc. Transl.} (1998).
\bibitem{Kon}
M.Kontsevich.
Vassiliev's knot invariants.
{\it Adv.in Sov.Math.} {\bf 16} (1993), pp.137-150.
\end{thebibliography}
\end{document}